\documentclass[12pt]{article}

\usepackage{amssymb}
\usepackage{epsfig}
\usepackage{amsfonts}
\usepackage{amsmath}
\usepackage{amssymb}

\newcommand*{\pd}[2]{\frac{\partial #1}{\partial #2}}

\begin{document}




\title{On a two-phase Hele-Shaw problem with a time-dependent gap and distributions of sinks and sources}
\author{T.V.~Savina, L.~Akinyemi, and A.~Savin}


\maketitle

\begin{abstract} 
A  two-phase Hele-Show problem with a time-dependent gap describes the evolution of the interface, which separates two fluids sandwiched between two plates. The fluids have different viscosities.  
  In addition to the change in the gap width of the Hele-Shaw cell, the interface   is driven by
  the presence of some special distributions of sinks and sources located in  both the interior and exterior domains. 
The effect of surface tension is neglected.
Using the Schwarz function approach,
we give examples of exact solutions when the interface belongs to  a certain family of algebraic curves and the curves do not form cusps. The family of curves are defined by the initial shape of the free boundary. 
\end{abstract}

Muskat problem, Generalized Hele-Shaw flow, Schwarz function, Mother body.


\maketitle


\section{Introduction}

Free boundary problems have been  a significant part of modern mathematics for more than a century, since the celebrated Stefan problem, which describes solidification, that is, an evolution of the moving front between liquid and solid phases. Free boundary problems also appear in fluid dynamics,  geometry, finance, and many other applications  (see \cite{chen} for a  detailed discussion).
Recently, they started to play an important role in modeling of biological processes
involving moving  fronts  of  populations  or  tumors \cite{fri2015}. 
These processes include cancer, biofilms, wound healing,  granulomas, and atherosclerosis \cite{fri2015}.
Biofilms are defined as communities of microorganisms, typically bacteria,
that are attached to a surface. The biofilms motivated  Friedman et al \cite{fri2014}
to consider a  two-phase
free boundary problem, where  one phase is an incompressible viscous fluid,
and  the other phase  is a mixture of two incompressible fluids, which
represent the viscous fluid and the polymeric network (with bacteria attached to it)
associated with a biofilm.
Free boundary problems are also  used in modeling of a tumor growth with one phase to be the tumor region, and the other phase to be the normal tissue surrounding the tumor \cite{fri2013}.

A Muskat problem is a free boundary problem related to the theory of flows in porous media \cite{mus}.
It describes an evolution of an interface between two immiscible fluids, `oil' and `water',  in a Hele-Shaw cell or in a porous medium.
Here we study a two-phase Hele-Shaw flow assuming that the upper plate uniformly moves up or down changing the gap width of a Hele-Shaw cell.
 Hele-Shaw free boundary problems have been extensively studied over the last century (see \cite{Vas2009}, \cite{Vas2015} and references therein). There are two classical formulations of the Hele-Shaw problems: the one-phase problem, when one of the fluids is assumed to be viscous while the other is effectively inviscid (the pressure there is constant), and the two-phase (or Muskat) problem. 
A statement of the problem with a time-dependent gap between the plates was mentioned in 
\cite{EEK} among other generalized Hele-Shaw flows. The one-phase (interior) version of this problem was considered in \cite{tian}, where  conditions of existence, uniqueness, and regularity of solutions were established under assumption 
that surface tension effects  on the free boundary are negligible;  some exact solutions were constructed as well.
An interior problem with a time-dependent gap and a non-zero surface tension was considered 
in \cite{jpA2015}, where asymptotic solutions were obtained for the case when initial shape of the droplet is a weakly distorted circle.
Note also that the mathematical formulation of the interior
problem with a time-dependent gap is similar to the problem of evaporation of a thin
film \cite{agam}. When the surface tension is negligible, the pressure in both formulations
can be obtained as a solution to the Poisson's equation in a bounded domain 
with homogeneous Dirichlet data on the free boundary.

Much less progress has been made for the Muskat problem.  Regarding the problem with a constant gap width, we should mention works \cite{howison2000}-\cite{contExact}. Specifically,
Howison \cite{howison2000} has obtained several simple solutions including the traveling-wave solutions and the stagnation point flow. In \cite{howison2000}, an idea of a method for solving some two-phase problems was  proposed and used to reappraise the Jacquard-S\'eguier solution \cite{JS}.
Global existence of solutions to some specific two-phase problems was considered in 
\cite{FT}-\cite{YT}.
Crowdy \cite{crowdy2006} presented an exact solution to the Muskat problem for the elliptical initial interface between two fluids of different viscosity. 
In \cite{crowdy2006}, it was shown that an  elliptical inclusion of one fluid remains elliptical when placed in a linear ambient flow of another fluid.
 In  \cite{contExact}, new exact solutions to the Muskat problem were constructed, extending the results obtained in \cite{crowdy2006}, to other types of inclusions.  
This paper is concerned with a two-phase Hele-Shaw problem with a variable gap width in the presence of sinks and sources.

Let $\Omega _2 (t) \subset {\mathbb R}^2$ with a boundary $\Gamma (t)$ at time $t$ be a 
simply-connected bounded domain occupied by a fluid with a constant viscosity $\nu _2$,   and let $\Omega _1 (t)$ be the region  ${\mathbb R}^2\setminus {\bar \Omega}_2(t)$ occupied by a different fluid of viscosity $\nu _1$.
To consider a two-phase Hele-Shaw flow forced by a time-dependent gap, we start with the Darcy's law 
\begin{equation}\label{1}
{\bf v}_j=-k_j\nabla p_j 
\quad\mbox{in}\quad \Omega _j (t), \qquad j=1,2,
\end{equation}
where ${\bf v}_j$ and $p_j$ are a two-dimensional gap-averaged velocity vector and a pressure of fluid $j$ respectively,
 $k_j=h^2(t)/12\nu _j$, and $h(t)$ is the gap width of the Hele-Shaw cell.
Equation (\ref{1}) is complemented by the volume conservation,   $A(t)h(t)=A(0)h(0)$ for any time $t$, where $A(t)$ and $A(0)$ are the areas of $\Omega _2 (t)$ and $\Omega _2(0)$ respectively.
The conservation of volume for a time-dependent gap may be written as a modification of the usual incompressibility condition 
$$
\nabla\cdot {\bf V_2} =0,
$$
where ${\bf V_2}=(u,v,w)$ is a three-dimensional velocity vector of the fluid occupying the domain $\Omega _2 (t)$. 
Indeed, the averaging of the three-dimensional incompressibility condition across the gap gives \cite{tian}:
$$
0=\int\limits _{0}^{h(t)} (u_x+v_y+w_z)dz/h(t)=u_x^{av}+v_y^{av}+(w(h(t))-w(0))/h(t)=u_x^{av}+v_y^{av}+\frac{\dot h(t)}{h(t)}.
$$
Here $z=0$ corresponds to the lower plate and $z=h(t)$ corresponds to the upper plate, and and $h(t)$ and $\dot h(t)$ are assumed to be small enough to avoid any inertial effects as well as to keep the large aspect ratio.
The latter implies \cite{tian}
\begin{equation}\label{2111}
\nabla\cdot {\bf v_2}=-\frac{\dot h(t)}{h(t)} \quad\mbox{in}\quad \Omega (t).
\end{equation}
Note that similar consideration may be applied to any finite part of the region $\Omega _1(t)$. 
Thus, equations (\ref{1}) and (\ref{2111}) suggest  to formulated the problem in terms of  
the pressure $p_j$ as a solution to  Poisson's equation,
\begin{equation}\label{01}
\Delta p_j=\frac{1}{k_j}\frac{\dot h(t)}{h(t)},
\end{equation}
 almost everywhere in the region $\Omega _j (t)$, satisfying boundary conditions
\begin{eqnarray}\label{2}
p_1(x,y,t)=p_2(x,y,t) \quad \mbox{on} \quad \Gamma (t),\\
-k_1\pd{p_1}{n}=-k_2\pd{p_2}{n} =v_n \quad \mbox{on} \quad \Gamma(t).\label{3}
\end{eqnarray}

We remark that when sinks and sources are present in $\Omega _j (t)$, equation (\ref{01}) has an additional term, $\Delta p_j=\frac{1}{k_j}\frac{\dot h(t)}{h(t)}+\mu_j$, describing the corresponding distribution.
Equation (\ref{2}) states the continuity of the pressure under the assumption of negligible
 surface tension.  
 Equation (\ref{3}) means that the normal velocity of the boundary itself coincides with the normal velocity of the fluid
 at the boundary.

The free boundary $\Gamma (t)$ moves due to a change of the gap width as well as the presence of  sinks and sources located in both regions. The supports of the sinks and sources,
specified in section \ref{prelim}, are either points or lines/curves. The presence of sinks and sources obviously changes the dynamics of the evolution of the interface between the fluids, which is shown for an elliptical interface in section \ref{sec:ex}.

For what follows, it is convenient to reformulate the problem in terms of harmonic functions $\tilde p_j$, where
\begin{equation}
p_j(x,y,t)=\tilde p_j(x,y,t)+\frac{1}{4k_j}\frac{\dot h(t)}{h(t)}(x^2+y^2).
\end{equation}
Then the problem \eqref{01}-\eqref{2} reduces to
\begin{equation}\label{01t}
\Delta \tilde p_j=\chi _j\mu _j \quad \mbox{in} \quad \Omega _j (t),
\end{equation}
where $\chi _j=0$ or $\chi _j=1$ in the absence or presence  of sinks and sources in
$\Omega _j (t)$ respectively,
\begin{eqnarray}\label{2t}
\tilde p_1(x,y,t)=\tilde p_2(x,y,t)+ \frac{k_1-k_2}{4k_1k_2}\,\frac{\dot h(t)}{h(t)}(x^2+y^2)
\quad \mbox{on} \quad \Gamma (t),\\
-k_1\pd{\tilde p_1}{n}
=-k_2\pd{\tilde p_2}{n} =
v_n +\frac{1}{4}\frac{\dot h(t)}{h(t)}\pd{}{n}(x^2+y^2)
\quad \mbox{on} \quad \Gamma(t).\label{3t}
\end{eqnarray}

The main difficulty of the two-phase problems is the fact that the pressure on the interface is unknown. However, if we assume that the free boundary remains within the family of curves, specified by the initial shape of the interface separating the fluids (which is feasible if the surface tension is negligible), the problem is drastically simplified. 

In this paper, using reformulation of the Muskat problem with the time-dependent gap in terms of the Schwarz function equation,  we describe a method of constructing exact solutions, and using this method we consider examples in the presence and in the absence of additional sinks and sources.

The structure of the paper is as follows.    In Section \ref{prelim} 
we describe the method
 of finding exact solutions. 
Examples of the exact solutions are given in Section \ref{sec:ex}, and concluding remarks are given in Section \ref{sec:concl}.

\section{The  method of finding exact solutions for a Muskat problem with a time-dependent gap}\label{prelim}

Consider a problem
\begin{equation}\label{01tm}
\Delta \tilde p_j=\chi _j\mu _j \quad \mbox{in} \quad \Omega _j (t),
\end{equation}
\begin{eqnarray}\label{2tm}
\tilde p_1(x,y,t)+\Psi _1(x,y,t)=\tilde p_2(x,y,t)+ \Psi _2(x,y,t)
\quad \mbox{on} \quad \Gamma (t),\\
-k_1\pd{\tilde p_1}{n}
=-k_2\pd{\tilde p_2}{n} =
v_n +\Phi (x,y,t)
\quad \mbox{on} \quad \Gamma(t).\label{3tm}
\end{eqnarray} 
In the case when
\begin{equation}\label{Psi}
\Psi _j=\frac{1}{4k_j}\,\frac{\dot h(t)}{h(t)}(x^2+y^2), \qquad j=1,2,
\end{equation}
\begin{equation}\label{Phi}
\Phi =\frac{1}{4}\,\frac{\dot h(t)}{h(t)}\pd{}{n}(x^2+y^2), \qquad j=1,2,
\end{equation}
the problem \eqref{01tm}-\eqref{3tm} coincides with \eqref{01t}-\eqref{3t}.

As stated before, the 
evolution of the interface separating the fluids is forced by the change in the gap width and the presence  of sinks and sources. In the absence of the surface tension, there is a possibility to control the interface by keeping 
$\Gamma(t)$ within a  family of curves defined by $\Gamma (0)$. 
For what follows, it is convenient to reformulate problem   (\ref{01tm})--(\ref{3tm}) in terms of the Schwarz function $S(z,t)$ of the curve $\Gamma (t)$ \cite{davis}--\cite{shapiro}.
This function for
 a real-analytic curve  $\Gamma :=\{ g(x,\,y,\,t)=0\}$  is defined as a solution 
to the equation
$g\left ((z+\bar z)/2, \, (z-\bar z)/2i,\, t \right ) =0$
with respect to $\bar z$.
This (regular) solution exists in some neighborhood $U_{\Gamma}$ of the curve $\Gamma$,
if the assumptions  of the implicit function theorem are satisfied \cite{davis}.
Note that if $g$ is a polynomial, then the Schwarz function is continuable into $\Omega _j$,
generally as a multiple-valued analytic function with a finite number of algebraic singularities (and poles).
In $U_{\Gamma}$, the normal velocity, $v_n$, of  $\Gamma (t)$ can be written in terms of the Schwarz function \cite{howison92},
$v_n=-i\dot S(z,t)/\sqrt{4\partial _z S(z,t)}$.

Let $\tau$ be an arclength along $\Gamma (t)$, $\psi _j$ be a stream function, and 
$W_j=\tilde p_j-i\psi _j$ be the complex potential, that is defined  on $\Gamma (t)$ and in $\Omega _j (t)\cap U_{\Gamma}$, $j=1,\,2$. Following \cite{cummings}-\cite{mcdonald},
taking into account the Cauchy-Riemann conditions in the $(n,\tau )$ coordinates, for the derivative of $W_j(z,t)$ with respect to $z$ on $\Gamma (t)$ we have
\begin{equation}\label{main1}
\partial _z {W_j}=\frac{\partial _{\tau}W_j}{\partial _{\tau}z}=
\frac{\partial _{\tau}\tilde p_j+i\partial _{n}\tilde p_j}{\partial _{\tau}z}=
\frac{\partial _{\tau}\tilde p_j-i(v_n+\Phi)/k_j}{\partial _{\tau}z}.
\end{equation}
Expressing $\partial _{\tau}z$ in terms of the Schwarz function, $\partial _{\tau}z=(\partial _z S(z,t))^{-1/2}$, we obtain 
\begin{equation}\label{main1m}
\partial _z {W_j}=\partial _{\tau}\tilde p_j\sqrt{\partial _z S}
-\frac{\dot {S}}{2k_j}-\frac{i\Phi}{k_j}\sqrt{\partial _z S}.
\end{equation}
Here $\partial _z {W}_j\equiv \pd{W_j}{z}$,  $\partial _{\tau} {}\equiv \pd{}{\tau}$.
Equation \eqref{2tm} implies that $\tilde p_1+\Psi_1=\tilde p_2+\Psi _2=f$ on $\Gamma (t)$,  where $f$ is an unknown function.
To keep $\Gamma (t)$ in a certain family  of curves defined by $\Gamma (0)$, for example, in a family of ellipses, we  assume that $f$ on $\Gamma (t)$ is a function of time only. 
This possibility is shown in Section 3, where specific examples are discussed.
In that case the problem is simplified drastically, and on $\Gamma (t)$ we have
\begin{equation}\label{main2m}
\partial _z {W_j}=
-\frac{\dot {S}}{2k_j}-\partial _z (\Psi _j(z,S(z,t))-\frac{i\Phi}{k_j}\sqrt{\partial _z S}\qquad j=1,2.
\end{equation}
For the special case when $\Psi _j$ and $\Phi$ are given by \eqref{Psi}, \eqref{Phi}, the last equation reduces to
\begin{equation}\label{main2mg}
\partial _z {W_j}=
-\frac{1}{2k_j}(\dot {S}+\frac{\dot h}{h} S)\qquad j=1,2.
\end{equation}

Remark that
each equation \eqref{main2mg} can be continued off of $\Gamma$ into the corresponding
$\Omega _j$, where $W_j$ is a multiple-valued analytic function. 
The equations \eqref{main2m} and  \eqref{main2mg}  imply that the singularities of $W_1$, $W_2$, and the singularities of the Schwarz function are linked. 
As such, the singularities of the Schwarz function play the crucial role in the construction of solutions in question.

To find the exact solutions,
suppose that at $t=0$ the interface is an algebraic curve, $\sum _{k=0}^n a_k(0) x^{k-n}y^n=0$, with the
Schwarz function $S(z,a_k^0)$. Assume  that during the course of evolution the Schwarz function of the interface $S(z,a_k(t))\equiv S(z,t)$ is such that
$S(z,a_k(0))=S(z,a_k^0)$, which leads us to the following six  
steps  method:
\newline
1) Compute $\dot S (z, t)$, locate its singularities, and define their type.
\newline
2) Using equations \eqref{main2mg} find  preliminary expressions for $\partial _z W_j$. 
\newline
3)  By putting restrictions on the coefficients $a_k(t)$ in the preliminary 
 expressions for $\partial _z W_j$
eliminate the terms involving undesirable singularities (if possible).
\newline
4)  Integrate \eqref{main2mg} with respect to $z$ in order to find   $W_j$ up to an arbitrary function of time.
\newline
5) Take the real part of $W_j$ in order to obtain $p_j$ up to an arbitrary function of time.
\newline
6) Evaluate the quantities $p_j$ on the interface to determine the independent of $z$ function of integration from the steps 3 and 4.
\newline
7) Locate the supports and compute the distributions of sinks and sources.

Before describing how to locate
the supports, 
we remark that  the distributions in  step 7 are related to the two-phase mother body 
\cite{contExact}.
The notion of  a mother body arises from the potential theory \cite{Gu1}-\cite{gardiner}
 and was adopted to the one-phase Hele-Shaw problem in \cite{external}.

As mentioned above, generally, the complex potentials $W_j$ are multiple-valued functions in $\Omega _j$.
For instance, if $\Gamma (t)$ is an algebraic curve, then the singularities of  $W_j$ are either poles or algebraic singularities.
To choose a branch of $W_j$, one has to introduce the cuts, $\gamma _j(t)$, that serve as 
supports for the distributions of sinks and sources, $\mu _j(t)$, $j=1,2$. 
Thus, each  cut  originates from 
an algebraic singularity  $z_a(t)$ of the potential $W_j$.
 The supports  consist of  those cuts  and/or points and  do not bound any two-dimensional subdomains in 
$\Omega _j(t)$, $j=1,2$. 
Each cut included in the support of  $\mu _j(t)$ is contained in the domain $\Omega _j (t)$, and the limiting values of
the pressure on each side of the cut are equal. 
The value of the density of sinks and sources located
 on the cut is equal to the jump of the normal derivative $\partial _n p_j$ of the pressure $p_j$.
In order for the total flux through the sinks and sources to be finite,
all of the singularities of the function $W_j$ must have no
more than the logarithmic growth.


The location of $z_a(t)$, as well as the directions of the cuts emanating from $z_a(t)$, are determined by  the Schwarz function via \eqref{main2mg}. In the examples considered below, the Schwarz function has the following two representations near its  singular points.
The first representation being  the square root (general position)
\begin{equation}\label{eq6}
S^g\left( z,t\right) =\xi^g\left( z,t\right)\,\sqrt{z-z_a(t)} +\zeta^g\left(
z,t\right),
\end{equation}
where $z_a(t)$ is a non-stationary singularity, that is $\dot z_a\ne 0$. 
The second being the reciprocal square root
\begin{equation}\label{eq6r}
S^r\left( z,t\right) =\frac{\xi^r\left( z,t\right)}{\sqrt{z-z_a(0)}} +\zeta^r\left(
z,t\right),
\end{equation}
where $z_a(0)$ is a stationary singularity, that is $\dot z_a= 0$.
Here $\xi^{g,r}\left( z,t\right) $ and $\zeta^{g,r}\left( z,t\right) $ are regular functions of $z$ in a
neighborhood of the point $z_a(t)$, and $\xi^{g,r}\left( z_a(t),t\right)\ne 0$.

By plugging \eqref{eq6} and \eqref{eq6r} into \eqref{main2mg}, in a small neighborhood of 
$z_a(t)$ we have
\begin{equation}\label{eq6w}
W_j^g\left( z,t\right) =\frac{1}{2k_j}\dot z_a\xi^g\left( z_a(t),t\right)\,\sqrt{z-z_a(t)} +\dots ,
\end{equation}
\begin{equation}\label{eq6wr}
W_j^r\left( z,t\right) =\frac{1}{k_j}C_0(t)\,\sqrt{z-z_a(0)} +\dots ,
\end{equation}
where the dots correspond to the smaller and regular terms that do not affect the computation of the directions of the cuts. The quantity $C_0(t)$ is defined by
$$
C_0(t)=\dot\xi^r\left( z_a(0),t\right )+\frac{\dot h(t)}{h(t)}\xi^r\left( z_a(0),t\right).
$$
Formulas \eqref{eq6w} and \eqref{eq6wr} along with the substitutions 
$z=z_a+\rho \exp{(i\varphi ^{g,r})}$ (with small $\rho$), imply that
\begin{equation}\label{eq6p}
p_j^g\left( z,t\right) =\frac{\sqrt{\rho}}{2k_j}\Re [\dot z_a\xi^g\left( z_a(t),t\right)\,
\exp{(\frac{i\varphi ^{g}}{2})}] +\dots ,
\end{equation}
\begin{equation}\label{eq6pr}
p_j^r\left( z,t\right) =-\frac{\sqrt{\rho}}{k_j}\Re [C_0(t)\exp{(\frac{i\varphi ^{r}}{2})}] +\dots .
\end{equation}
Computing the zero level of a variation of $p_j$ along a small loop surrounding the singular point, we finally 
obtain the following directions of the cuts:
for the general position
\begin{equation}\label{dir}
\varphi ^g=\pi-2(\arg [\xi^g\left( z_a(t),t\right)]+ \arg [\dot z_a]) +2\pi k,\quad k=0,\pm1, \pm 2 ... .
\end{equation} 
and  for the reciprocal square root
\begin{equation}\label{dirr}
\varphi ^r=\pi-2\arg [C_0(t)] +2\pi k,\quad k=0,\pm1, \pm 2 ... .
\end{equation}

In the next section, we use the described method to construct  exact solutions to the Muskat problem. In the considered examples, the evolution of the interface is driven by the change in the  gap width of the Hele-Shaw cell. The examples include
the elliptical shape with and without sinks and sources in the finite domain as well as the Cassini's oval in the presence of sinks and sources.

\section{ Examples of specific initial interfaces \label{sec:ex}}

\subsection{Circle } 

To illustrate the method, we start with the simplest example for which the solution is known. Suppose that the initial shape of the interface is a circle with the equation $x^2+y^2=a^2(0)$, and during the evolution the boundary remains circular, 
 $x^2+y^2=a^2(t)$. The corresponding Schwarz function is $S=a^2(t)/z$. 
Taking into account the volume conservation, equation \eqref{main2mg} in this case reads as
$\partial _z W_j=0$, which implies that $\tilde p_j$ is a function depending on $t$ only,
\begin{equation}
\tilde p_j=-\frac{a_0^2h_0\dot h}{4k_jh^2}+f(t),
\end{equation}
therefore,
\begin{equation}\label{Pcircle}
p_j(x,y,t)=\frac{1}{4k_j}\frac{\dot h(t)}{h(t)}\Bigl ( x^2+y^2-\frac{a_0^2h_0}{h(t)}\Bigr )+f(t)
\end{equation}
and $a(t)=a_0\sqrt{h_0/h(t)}$.

\subsection{Ellipse}

Consider a two-phase problem with an elliptical interface,
$\Gamma (0)=\left\{ \frac{x^2}{a(0)^2}+\frac{y^2}{b(0)^2}=1\right\}$,
  where $a(0)$ and $b(0)$ are given and $a(0)>b(0)$.
The Schwarz function of an elliptical interface with  semi-axes $a(t)$ and $b(t)$ is
$$
S\left( z,t\right)
=\Bigl ( \bigl ( a(t)^2+b(t)^2 \bigr ) z-2a(t)b(t)\sqrt{z^2-d(t)^2}\Bigr )/d(t)^2,
$$
where $d(t)=\sqrt{a(t)^2-b(t)^2}\,$ is the half of the inter-focal distance.
Assuming that 
the interface remains elliptical during the course of the evolution, 
we use equation \eqref{main2mg}
\begin{align}\notag
  \partial_zW_j=-\frac{1}{2k_j}\left(\partial_t{S}+\frac{\dot{h}}{h}{S}\right).
\end{align}
Due to the volume conservation of the fluid occupying $\Omega _2(t)$, the product of functions $a(t)$ and $b(t)$ is linked to the gap width, $h(t)$, via the equation $h(t)=a_0b_0h_0/(a(t)b(t)),$ where $a_0=a(0)$, $b_0=b(0),$ and $h_0=h(0)$. Therefore, $\dot{h}(t)/{h(t)}=-\partial_t(ab)/(ab)$, and the equation \eqref{main2mg}  could be rewritten as 
\begin{align}\label{27.91}
   \partial_zW_j=-\frac{1}{2k_j}\left(\partial_t{S}-\frac{\partial_t(ab)}{ab}{S}\right),
\end{align}
which results in
\begin{eqnarray}\label{wzellipse}
\partial _z{W_j}=-\frac{z}{2k_j}\Bigl \{\pd{}{t}\Bigl (  \frac{a^2+b^2}{d^2}   \Bigr )
-\frac{(a^2+b^2)}{a\,b\, d^2}\pd{}{t}(ab)\Bigr \}
 \notag \\
-\frac{(2z^2-d^2)}{ \sqrt{z^2-d^2}}\, \frac{ab}{2k_jd^4}\pd{}{t}\Bigl (  d^2   \Bigr )
\end{eqnarray}
and
\begin{eqnarray}\label{wzellipse1}
W_j=-\frac{z^2}{4k_j}\Bigl \{\pd{}{t}\Bigl (  \frac{a^2+b^2}{d^2}   \Bigr )
-\frac{(a^2+b^2)}{a\,b\, d^2}\pd{}{t}(ab)\Bigr \}
 \notag \\
-\frac{a\,b\,z}{2k_jd^4}\sqrt{z^2-d^2}\,\pd{}{t}(d^2) 
+C_j(t),\label{well}
\end{eqnarray}
where $C_j(t)$ is an arbitrary function of time.
\begin{figure}
\hskip1.5cm
\begin{center}
\includegraphics[width=0.3\textwidth]{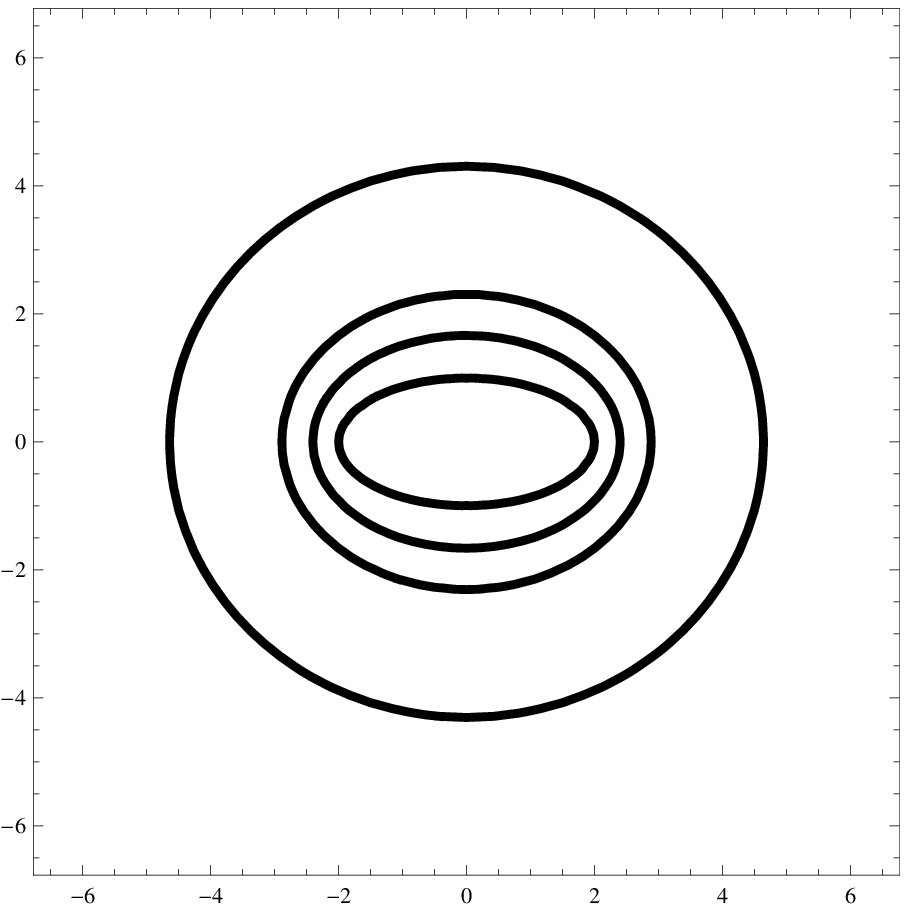}{a}
\includegraphics[width=0.3\textwidth]{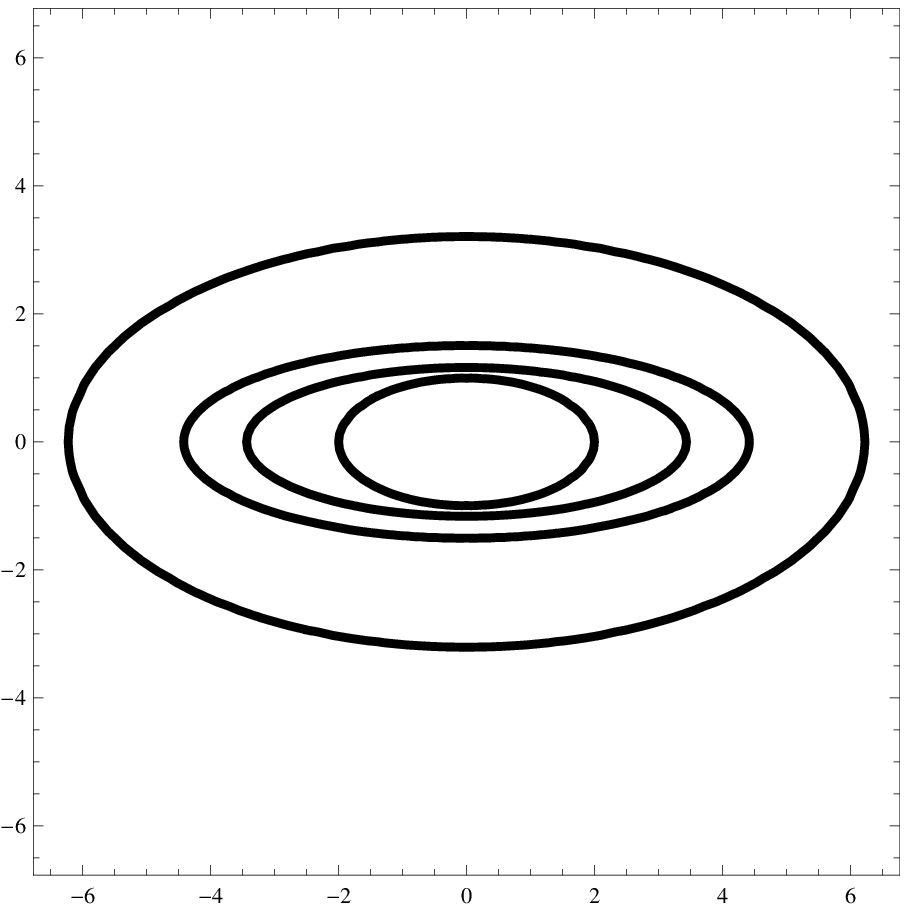}{b}
\includegraphics[width=0.3\textwidth]{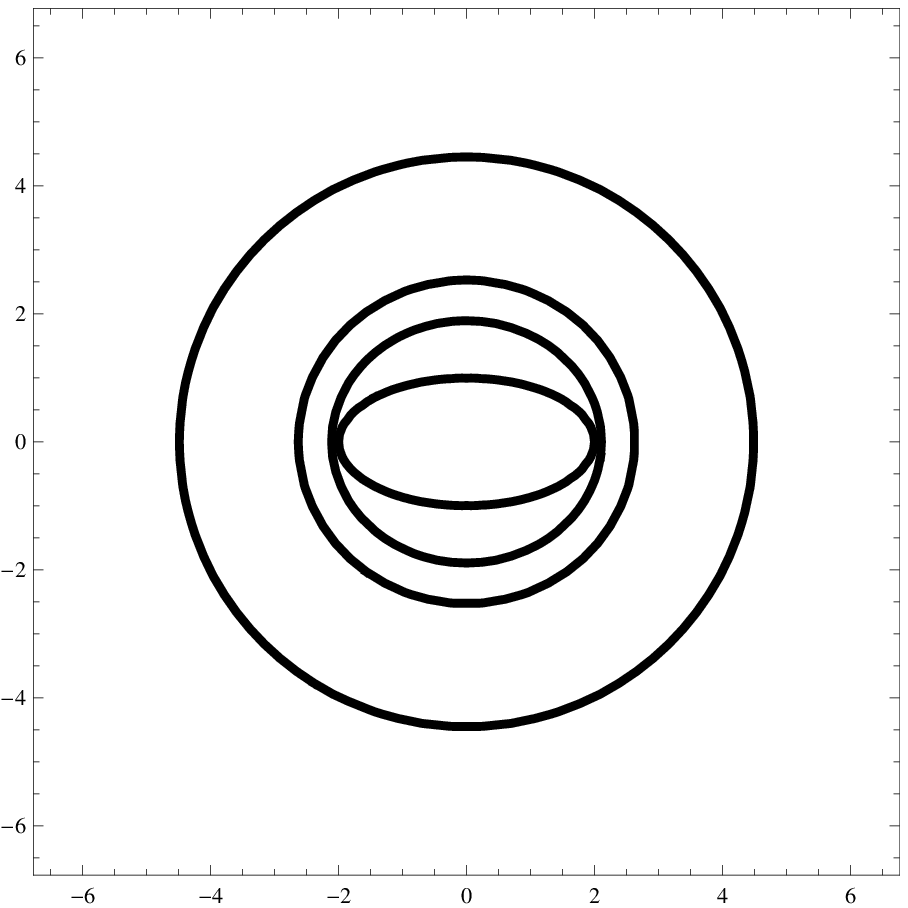}{c}
\caption{Squeezing of an ellipse:  $a_0=2$, $b_0=1$, $h_0=0.1$, $h(t)=h_0-t$; $t=0$, $t=0.05$, $t=0.07$, $t=0.09$:
(a) $d^2=const$, (b) $d^2(t)=d_0^2\exp{(25t)}$,
(c) $d^2(t)=d_0^2\exp{(-25t)}$. }
\label{fgElli}
\end{center}
\end{figure}

{\it (a) Evolution with constant inter-focal distance}.

To obtain an exact solution in the absence of sinks and sources in the finite part of the plane, we  set $d(t)=d(0)$. Then, the second term in the formula (\ref{wzellipse1}) vanishes, which implies the following expression for the pressure 
\begin{equation}
\tilde p_j =\Re [W_j]=\frac{1}{4k_j}\Bigl ( (x^2-y^2)\frac{\dot a\, d_0^2}{a(a^2-d_0^2)}+2\dot a a  
 \Bigr )+f(t),
\end{equation}
therefore,
\begin{equation}
p_j =\frac{\dot a}{2k_j a (a^2-d^2_0)}\Bigl ( d^2_0 x^2-a^2(x^2+y^2)+a^2(a^2-d^2_0) 
 \Bigr )+f(t)
\end{equation}
is the solution to the problem \eqref{01}-\eqref{3}.
Note that when $d_0=0$, this formula coincides with formula \eqref{Pcircle} related to the circular interface.

Hence, $\Gamma (t)$ is a family of co-focal ellipses, 
$$
\frac{x^2}{a^2(t)}+\frac{y^2}{b^2(t)}= 1,
$$
controlled by one of the functions $a(t)$, $b(t)$  or $h(t)$. If $h(t)$ is given, then
\begin{eqnarray}
& a^2(t)=\frac{1}{2}\Bigl ( a_0^2-b_0^2+\sqrt{(a_0^2-b_0^2)^2+4a_0^2b_0^2h_0^2/h^2(t)}    \Bigr ),\\
& b^2(t)=\frac{1}{2}\Bigl ( b_0^2-a_0^2+\sqrt{(a_0^2-b_0^2)^2+4a_0^2b_0^2h_0^2/h^2(t)}    \Bigr ).
\end{eqnarray}
An example of such an evolution with  a linear function $h(t)$ is shown in Fig.~1(a).

{\it (b) Evolution with variable inter-focal distance}.

If we admit solutions with variable inter-focal distance by keeping all terms in 
(\ref{wzellipse1}), we must allow, in addition to the gap change, some sinks/sources located in $\Omega _2$.
In that case, the pressure is
\begin{eqnarray}
\tilde p_j=-\frac{(x^2-y^2)}{4k_j}\Bigl \{\pd{}{t}\Bigl (  \frac{a^2+b^2}{d^2}   \Bigr )
-\frac{(a^2+b^2)}{a\,b\, d^2}\pd{}{t}(ab)\Bigr \}
 \notag \\
-\frac{a\,b}{2k_jd^4}\,\pd{}{t}(d^2) \frac{x\,(\alpha ^2-y^2)}{\alpha}
-\frac{ab\,(\dot ab-a\dot b)}{2k_jd^2}
+f(t),\label{well}
\end{eqnarray}
where 
$$
\alpha ^2=\Bigl ( x^2-y^2-d^2 +\sqrt{(x^2-y^2-d^2)^2+4x^2y^2}  \Bigr ) /2,
$$
therefore, making
\begin{eqnarray}
 p_j=-\frac{(x^2-y^2)}{4k_j}\Bigl \{\pd{}{t}\Bigl (  \frac{a^2+b^2}{d^2}   \Bigr )
-\frac{(a^2+b^2)}{a\,b\, d^2}\pd{}{t}(ab)\Bigr \}
-\frac{ab\,(\dot ab-a\dot b)}{2k_jd^2}
\notag \\
-\frac{a\,b}{2k_jd^4}\,\pd{}{t}(d^2) \frac{x\,(\alpha ^2-y^2)}{\alpha}
-\frac{\partial _t(ab)}{4k_j ab}(x^2+y^2)
+f(t).\label{well1}
\end{eqnarray}
Equation (\ref{wzellipse1}) implies that there are two singular points in the interior domain $\Omega _2$, $z=\pm d$. The Schwarz function near those points has the square root representation \eqref{eq6} with 
$$
\xi ^g=-\frac{2ab}{d^2}\sqrt{z\pm d}.
$$
The direction of the cut at each point is defined by formula \eqref{dir}, 
which implies that 
at the point $z_a=d$, the angle is
$\varphi ^g =\pi +2\pi k$  and at the point $z_a=-d$,  the angle is $\varphi ^g=2\pi k$, $k=0,\pm 1,\pm 2,\dots$. Thus, the cut $\gamma _2(t)$ is located along the inter-focal segment $[-d,d]$.
The density of the distribution of sinks and sources along that segment is given by the formula
$$
\mu _2= \frac{ab\,\partial _t (d^2)}{k_2d^4}\,\,\frac{(2x^2-d^2)}{\sqrt{d^2-x^2}}.
$$
 Such a density changes its sign along the inter-focal segment, so its presence does not affect the area of the ellipse, 
 $\dot A=\int _{-d}^{d} k_2 \mu _2(x,t)\, dx=0$.
Fig.~1 shows how the sinks and sources change the evolution of the interface with increasing (see Fig.~1 (b)) and decreasing (see Fig.~1 (c)) inter-focal distances. 

\subsection{The  Cassini's oval}

\begin{figure}
\hskip1.5cm
\begin{center}
\includegraphics[width=0.3\textwidth]{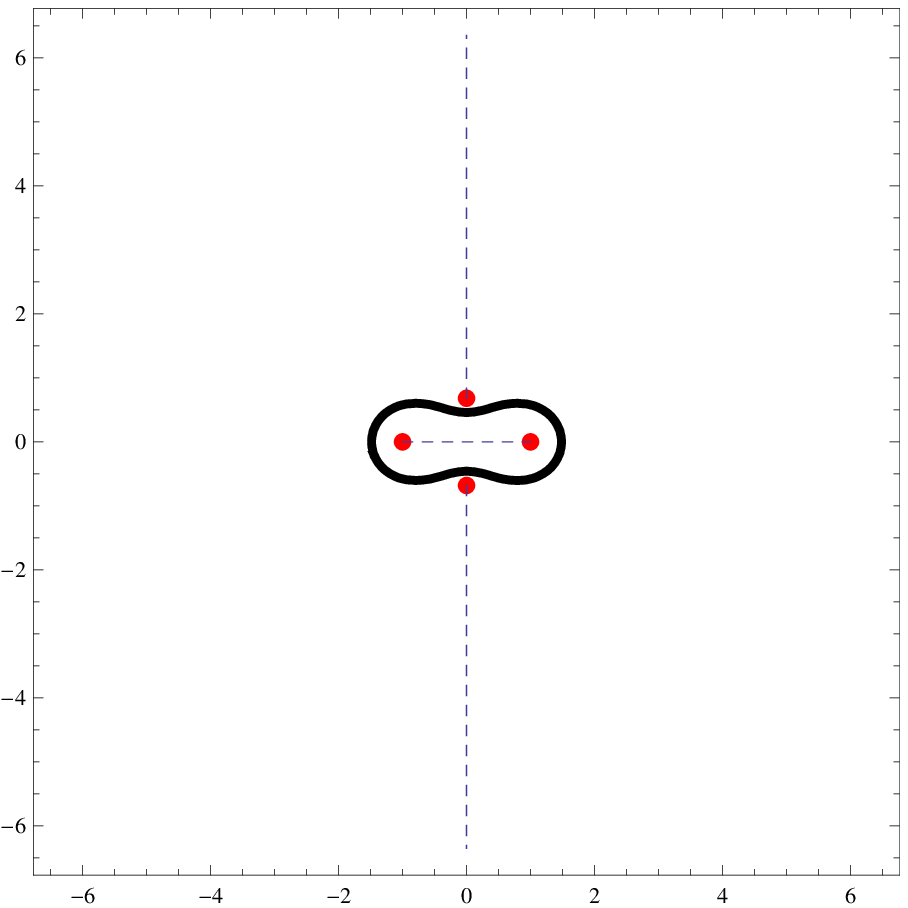}{a}
\includegraphics[width=0.3\textwidth]{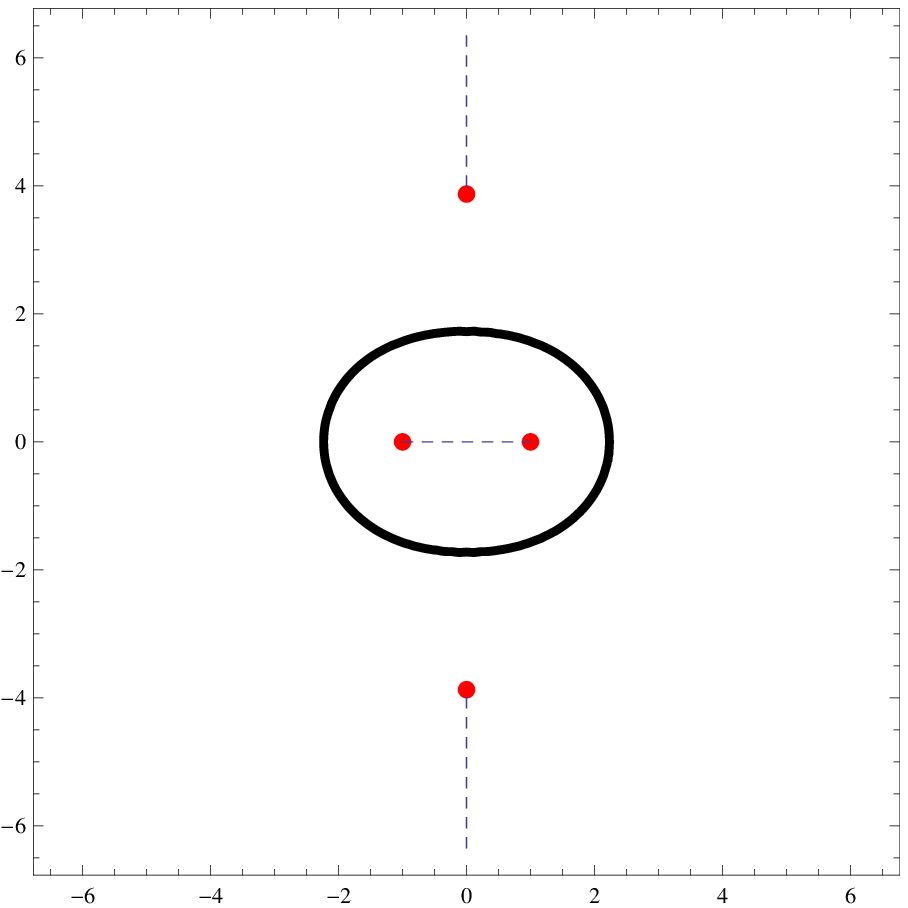}{b}
\caption{Squeezing of the Cassini's ovals  
for $b(t)=b_0=1$, $a_0=1.1$, $h_0=0.1$,
$h(t)=h_0-t$: 
(a) $t=0$, 
 (b) $t=0.05$. 
} 
\label{fgCass}
\end{center}
\end{figure}

Similar to the previous examples, assume that $\Gamma (t)$   remains in the specific family of curves,
the Cassini's ovals, given by the equation
$$
\left( x^2+y^2\right) ^2-2b(t)^2\left( x^2-y^2\right) =a(t)^4-b(t)^4,
$$
where $a(t)$ and $b(t)$ are unknown positive functions of time.
This curve consists  of one
closed curve, if $a(t)>b(t)$ (see Fig. \ref{fgCass}), and  two closed curves otherwise. Assume that at $t=0$ $a(0)>b(0)$. 
The Schwarz function of  Cassini's oval,
$$
S\left( z,t \right) =\sqrt{b^2z^2+a^4-b^4}\, /\sqrt{z^2-b^2},
$$
 has two singularities in $\Omega _1(t)$, $z=\pm i\sqrt{(a^4-b^4)/b^2}$, and two singularities  in $\Omega _2(t)$, $z=\pm b$.
The corresponding complex velocities have singularities at the same points,
\begin{equation}\label{I3}
\partial _z W_j=-\frac{1}{2k_j}\Bigl ( 
\frac{B_1z^2+B_2}{\sqrt{(b^2z^2+a^4-b^4)(z^2-b^2)}}+
\frac{b\dot b \sqrt{b^2z^2+a^4-b^4}}{\sqrt{(z^2-b^2)^3}}
 \Bigr ).
\end{equation}
Here
$$
B_1=b\dot b +b^2\dot h /h,\qquad
B_2 =2a^3\dot a-2b^3\dot b + (a^4-b^4)\dot h/h,
$$
and $\dot h/h=-\dot A/A$ due to volume conservation.

The area of   Cassini's  oval can be computed in polar coordinates, 
$A=a^2E(\pi,\frac{b^2}{a^2})=2a^2E(\frac{b^2}{a^2})$, where 
$E(\phi, k)=\int\limits _0^\phi\sqrt{1-k^2\sin ^2 t }\,dt$
and $E(k)=E(\pi/2,k)$, resulting in
\begin{equation}
\frac{\dot A}{A}=\frac{2\dot a}{a}+\frac{\partial _t E(\pi,\frac{b^2}{a^2})}{E(\pi,\frac{b^2}{a^2})}.
\end{equation}
Taking into account (\cite{prudnikov}, p. 772),
$$
\pd{E(\phi,k)}{k}=\frac{1}{k}\Bigl ( E(\phi ,k)-F(\phi, k)\Bigr ),
$$
where 
\begin{equation}\label{ellipticF}
F(\phi, k)=\int\limits _0^\phi\frac{1}{\sqrt{1-k^2\sin ^2 t }}\,dt,
\end{equation}
$F(\pi/2,k)=K(k)$, and
$
\partial _t E(\pi,\frac{b^2}{a^2})=\Bigl (E(\pi,\frac{b^2}{a^2})-F(\pi,\frac{b^2}{a^2})\Bigr )
\frac{2a\dot b-2b\dot a}{ab},
$
we have
\begin{equation}
B_1(t)=\frac{b}{aE(\pi,\frac{b^2}{a^2})}
\Bigl (-a\dot b E(\pi,\frac{b^2}{a^2})+2(a\dot b-\dot a b)F(\pi,\frac{b^2}{a^2})\Bigr ),
\end{equation}
\begin{equation}
B_2(t)=\frac{2(\dot ab-a\dot b)}{abE(\pi,\frac{b^2}{a^2})}
\Bigl (a^4 E(\pi,\frac{b^2}{a^2})-(a^4- b^4)F(\pi,\frac{b^2}{a^2})\Bigr ),
\end{equation}
and
\begin{equation}
 W_j=-\frac{1}{2k_j}\Bigl ( 
B_1 I_1+B_2I_2+
b\dot b I_3
 \Bigr ).
\end{equation}
Here
$$
I_1 =\frac{b^2}{a^2}F \Bigl (\cos ^{-1} \bigl ( \frac{b}{z}\bigr ),\, \frac{\sqrt{a^4-b^4}}{a^2}\, \Bigr )
-\frac{a^2}{b^2}E \Bigl (\cos ^{-1} \bigl ( \frac{b}{z}\bigr ),\, \frac{\sqrt{a^4-b^4}}{a^2}\, \Bigr )
+\frac{\sqrt{(z^2b^2+a^4-b^4)(z^2-b^2)}}{zb^2},
$$
\begin{eqnarray}\notag
& I_2 =\frac{1}{a^2}F \Bigl (\cos ^{-1} \bigl ( \frac{b}{z}\bigr ),\, \frac{\sqrt{a^4-b^4}}{a^2}\, \Bigr )=\frac{1}{a^2}\int\limits _0 ^{\sqrt{1-b^2/z^2}}
\frac{dt}{\sqrt{1-\frac{a^4-b^4}{a^4}t^2}\,\sqrt{1-t^2}}= \\
& \frac{1}{a^2}\int\limits _0 ^{\cos ^{-1}  (b/z )}\frac{dt}{\sqrt{1-\frac{a^4-b^4}{a^4}\sin ^2 t}}, \notag
\end{eqnarray}
and the integral $I_3$ corresponds to the last term in \eqref{I3}.
To ensure that the singularities of the complex potential have no more than the logarithmic type, we eliminate this term by setting $\dot b$ to zero.  Thus, we have
$$
\dot S\left( z\right) =\frac{2a^3\dot a}{\sqrt{b^2z^2+a^4-b^4}\sqrt{z^2-b^2}},
$$
and the equation \eqref{main2mg} implies
\begin{eqnarray}\label{mainCas}
& {W_j} =-\frac{a\dot a}{k_jE(\frac{b^2}{a^2})}\Bigl[\Bigl (E(\frac{b^2}{a^2})-K(\frac{b^2}{a^2})\Bigr )
\,F \Bigl (\xi,\, \frac{\sqrt{a^4-b^4}}{a^2}\,\Bigr )\\
& +K(\frac{b^2}{a^2})\,E \Bigl (\xi,\, \frac{\sqrt{a^4-b^4}}{a^2}\,\Bigr )
-\frac{K(\frac{b^2}{a^2})\,\sqrt{(z^2b^2+a^4-b^4)(z^2-b^2)}}{a^2z}\Bigr ]
+C(t), \notag
\end{eqnarray}
 where $\xi=\cos ^{-1}  \bigl (\frac{b}{z} \bigr)$ and $F(\alpha\,,\beta)$ is the incomplete elliptic integral of the first kind \eqref{ellipticF},
$$
F \Bigl (\cos ^{-1} \bigl ( \frac{b}{z}\bigr ),\, \frac{\sqrt{a^4-b^4}}{a^2}\, \Bigr )=\int\limits _0 ^{\sqrt{1-b^2/z^2}}
\frac{dt}{\sqrt{1-\frac{a^4-b^4}{a^4}t^2}\,\sqrt{1-t^2}}=
\int\limits _0 ^{\cos ^{-1}  (b/z )}\frac{dt}{\sqrt{1-\frac{a^4-b^4}{a^4}\sin ^2 t}}.
$$
Since $p_j=\Re \,[W_j]$, we need to compute the real parts for each term in \eqref{mainCas}.
Using the property $\overline {F(\alpha\,,\beta)}=F(\overline{\alpha}\,,\beta)$ and the summation formula for the elliptic integrals of the first kind \cite{BE},  we have
$$
\frac{1}{2}
\Bigl [F \Bigl (\xi,\, \frac{\sqrt{a^4-b^4}}{a^2}\,\Bigr )+\overline{F \Bigl (\xi,\, \frac{\sqrt{a^4-b^4}}{a^2}}\,\Bigr )\Bigr]=\frac{1}{2}F \Bigl (\alpha
,\, \frac{\sqrt{a^4-b^4}}{a^2}\,\Bigr ),
$$
where
\begin{equation}\label{alpha1}
\alpha=\sin ^{-1}\frac{
\cos \overline{\xi}\sin \xi \sqrt{1-\frac{{a^4-b^4}}{a^4}\sin ^2\overline{\xi}}+
\cos \xi\sin \overline{\xi} \sqrt{1-\frac{{a^4-b^4}}{a^4}\sin ^2{\xi}}
}
{   1-\frac{{a^4-b^4}}{a^4}\sin ^2 \xi \sin ^2 \overline{\xi}   }
\end{equation}
or
\begin{equation}\label{alpha3}
\alpha=\sin ^{-1}\frac{
a^2z\sqrt{z^2-b^2}\sqrt{b^2\bar z^2+a^4-b^4}+
a^2\bar z \sqrt{\bar z ^2 -b^2} \sqrt{b^2z^2+a^4-b^4}
}
{   b^2z^2\bar z ^2+(a^4-b^4)(z^2+\bar z ^2-b^2)  }.
\end{equation}
Similarly, using the property $\overline {E(\alpha\,,\beta)}=E(\overline{\alpha}\,,\beta)$ and the summation formula for the elliptic integrals of the second kind \cite{BE},  we have
$$
\frac{1}{2}
\Bigl [E \Bigl (\xi,\, \frac{\sqrt{a^4-b^4}}{a^2}\,\Bigr )+\overline{E \Bigl (\xi,\, \frac{\sqrt{a^4-b^4}}{a^2}}\,\Bigr )\Bigr]=\frac{1}{2}E \Bigl (\alpha
,\, \frac{\sqrt{a^4-b^4}}{a^2}\,\Bigr )+
\frac{(a^4-b^4)\sqrt{(z^2-b^2)(\bar z^2-b^2)}}{2a^4z\bar z}\sin\alpha    .
$$
Consequently, the pressure is determined by
\begin{eqnarray}\label{mainCasP}
{\tilde p_j} = &  -\frac{a\dot a}{2k_jE(\frac{b^2}{a^2})}\Bigl[\Bigl (E(\frac{b^2}{a^2})-K(\frac{b^2}{a^2})\Bigr )
\,F \Bigl (\alpha,\, \frac{\sqrt{a^4-b^4}}{a^2}\,\Bigr )
 +K(\frac{b^2}{a^2})\,E \Bigl (\alpha,\, \frac{\sqrt{a^4-b^4}}{a^2}\,\Bigr )\\
& +K(\frac{b^2}{a^2})\,\frac{(a^4-b^4)\sqrt{(z^2-b^2)(\bar z ^2-b^2)}}{a^4z\bar z}
-\frac{2K(\frac{b^2}{a^2})}{a^2}\,\Re\{\frac{\sqrt{(z^2b^2+a^4-b^4)(z^2-b^2)}}{z}\}\Bigr ]
+C_j(t). \notag
\end{eqnarray}
Here 
$$
\Re\{\frac{\sqrt{(z^2b^2+a^4-b^4)(z^2-b^2)}}{z}\}=
\frac{x(\alpha _1^2\alpha _2 ^2-x^2y^2b^2+y^2(\alpha _1^2b^2+\alpha _2^2))}
{(x^2+y^2)\,\alpha _1\alpha _2},
$$
where
$$
\alpha _1 ^2 = (x^2-y^2-b^2+\sqrt{(x^2-y^2-b^2)^2+4x^2y^2}\,)/2
$$
and
$$
\alpha _2 ^2 = ((x^2-y^2)b^2+a^4-b^4+\sqrt{((x^2-y^2)b^2+a^4-b^4)^2+4x^2y^2b^2}\,)/2.
$$
Taking into account the boundary condition to determine $C_j(t)$, we have
 \begin{eqnarray}\label{mainCasP1}
&  {\tilde p_j} =  -\frac{a\dot a}{2k_jE(\frac{b^2}{a^2})}\Bigl[\Bigl (E(\frac{b^2}{a^2})-K(\frac{b^2}{a^2})\Bigr )
\,F \Bigl (\alpha,\, \frac{\sqrt{a^4-b^4}}{a^2}\,\Bigr )   \notag
 +K(\frac{b^2}{a^2})\,E \Bigl (\alpha,\, \frac{\sqrt{a^4-b^4}}{a^2}\,\Bigr )\\
& + K(\frac{b^2}{a^2})\,\frac{(a^4-b^4)\sqrt{(z^2-b^2)(\bar z ^2-b^2)}}{a^4z\bar z}\notag
-\frac{2K(\frac{b^2}{a^2})}{a^2}\,\Re\{\frac{\sqrt{(z^2b^2+a^4-b^4)(z^2-b^2)}}{z}\}\\
& -\Bigl (E(\frac{b^2}{a^2})-K(\frac{b^2}{a^2})\Bigr )
\,K \Bigl (\frac{\sqrt{a^4-b^4}}{a^2}\,\Bigr )-
K(\frac{b^2}{a^2})\,E \Bigl (\frac{\sqrt{a^4-b^4}}{a^2}\,\Bigr )
\Bigr ]
+f(t) 
\end{eqnarray}
or
\begin{eqnarray}\label{mainCasP2}
&  {\tilde p_j} =  -\frac{a\dot a}{2k_jE(\frac{b^2}{a^2})}\Bigl[\Bigl (E(\frac{b^2}{a^2})-K(\frac{b^2}{a^2})\Bigr )
\,F \Bigl (\alpha,\, \frac{\sqrt{a^4-b^4}}{a^2}\,\Bigr )   \notag
 +K(\frac{b^2}{a^2})\,E \Bigl (\alpha,\, \frac{\sqrt{a^4-b^4}}{a^2}\,\Bigr )\\
& + K(\frac{b^2}{a^2})\,\frac{(a^4-b^4)\sqrt{(x^2+y^2)^2-2b^2(x^2-y^2)+b^4}}{a^4(x^2+y^2)}\\ \notag
& -\frac{2K(\frac{b^2}{a^2})}{a^2}\,\frac{x(\alpha _1^2\alpha _2 ^2-x^2y^2b^2+y^2(\alpha _1^2b^2+\alpha _2^2))}
{(x^2+y^2)\,\alpha _1\alpha _2}
 -\frac{\pi}{2}\Bigr ]
+f(t). \notag
\end{eqnarray}
Thereby,
$$
p_j=\tilde p_j- \frac{\dot a K(\frac{b^2}{a^2})}{2k_j aE(\frac{b^2}{a^2})}(x^2+y^2).
$$

To find the location of sinks and sources in the interior domain $\Omega _2$, note that the Schwarz function near its singular points $z=\pm b$ has the reciprocal square root representation  \eqref{eq6r} with
$\xi^r(z,t)=\sqrt{b^2z^2+a^4-b^4}/\sqrt{z\pm b}$. Formula \eqref{dirr} implies that 
$\varphi ^r(b)=\pi$ and $\varphi ^r(-b)=0$. This results (taking into account the symmetry of the problem) in the segment  $x\in [-b,b]$ as a location of sinks and sources.  
The corresponding density is
$$
\mu _2 =\frac{B_1x^2+B_2}{k_2\sqrt{(b^2x^2+a^4-b^4)(b^2-x^2)}}.
$$
Note that $\int\limits _{-b}^{b}\mu _2 (x)\, dx=0$, which  is consistent with the volume conservation.

To determine the location of the sinks and sources in  domain $\Omega _1$, we start with
singular points $z_a(t)= \pm i\sqrt{(a^4-b^4)}\,/b$.
The Schwarz function near these points has the square root representation \eqref{eq6}, and 
the directions of the cuts are defined by  formula \eqref{dir}. 

In the neighborhood of the point $z_a(t)= i\sqrt{(a^4-b^4)}\,/b$, we have $\arg [\dot z _a]=\pi /2 +2\pi k$ and
$\arg  [\xi ^g\left( z_a(t),t\right)]=-\pi /4+\pi  k$. Thus, according to  \eqref{dir}  the direction of the cut is $\varphi ^g=\pi /2+2\pi k$, $k=0,\pm 1,\pm2, \dots$.

Similarly, at the point $z_a(t)= -i\sqrt{(a^4-b^4)}\,/b$, $\arg [\dot z _a]=-\pi /2 +2\pi k$,
$\arg  [\xi ^g\left( z_a(t),t\right)]=-3\pi /4+\pi k$. Therefore, the direction of the cut is $\varphi ^g=-\pi /2+2\pi k$.

Taking into consideration symmetry with respect to the $x$-axis, we conclude that the 
 support  of $\mu _1$ consists of  two rays starting at the branch points and going to infinity (see the dashed lines in Fig.~2).  The density of sinks and sources is defined by 
$$
\mu _1 =\frac{B_1y^2-B_2}{k_1\sqrt{(b^2y^2-a^4+b^4)(b^2+y^2)}}.
$$
The evolution of the oval is controlled by a single function $h(t)$, where $b$ is constant and the parameter $a(t)$ is defined by the equation:
$$
\frac{\dot h}{h}=-\frac{\dot a}{a}\frac{K(b^2/a^2)}{E(b^2/a^2)}.
$$ 
Fig.~2 shows the evolution of the Cassini's oval under squeezing with $h(t)=h_0-t$ at $t=0$ (see Fig.~2 a) and $t=0.05$ (see Fig.~2 b). The dots correspond to the singular points $z_a$, the dashed lines correspond to the cuts.

\section{Concluding remarks}\label{sec:concl}

We have studied  a Muskat problem with a negligible surface tension and a gap width
dependent on time. This study extended the results reported in \cite{tian}, \cite{jpA2015}, and \cite{contExact}.
We suggested a method of finding exact solutions and applied it to find new exact solutions for initial elliptical shape and  Cassini's oval.   The idea of the  method was to keep the interface within a certain family of  curves defined by its initial shape.

For the elliptical shape,  we found two types of solutions: without sinks and sources in the interior domain, and with the presence of a special distribution of sinks and sources along the inter-focal distance. In the former solution, the  inter-focal distance remains constant,
while in the latter,  it changes. 

For the Cassini's oval, we found a solution to the problem when both a gap change and special distributions of sinks and sources in both the interior and exterior domains are present.

Our mathematical model included an assumption that the volume of the bounded domain 
$\Omega _2(t)$ is conserved. To show other conserved quantities, we follow Richardson \cite{rich}, \cite{EEK}, \cite{EE} deriving the moment dynamics equation,
\begin{equation}\label{moment}
\frac{d}{dt}\Bigl [h(t)\int\limits _{\Omega _2(t)} u(x,y) dxdy
\Bigr ]=-\chi _2k_2(t)h(t)\int\limits _{\gamma _2(t)} u(s)\mu _2 (s,t) \, ds, 
\end{equation}
where $u(x,y)$ is a harmonic function in a domain $\Omega \supset \Omega _2(t)$.
The latter follows from the chain of equalities:
$$
\frac{d}{dt}\Bigl [\int\limits _{\Omega _2(t)} u\, dxdy
\Bigr ]=\int\limits _{\Gamma(t)} u\,v_n ds
=-k_2\int\limits _{\Gamma(t)} u\,\pd{p_2}{n} d\tau
$$ 
$$
=-\frac{\dot h}{h}\int\limits _{\Omega _2(t)} u\, dxdy
-\chi _2k_2h\int\limits _{\gamma _2(t)} u\, \mu_2 ds-k_2f(t) \int\limits _{\Gamma(t)} 
\pd{u}{n} d\tau .
$$
By setting $f(t)$ to zero and rearranging the terms, we arrive at \eqref{moment}.

Equation \eqref{moment} implies that in the absence of sinks and sources, $\chi_2=0$, the quantity
$h\int _{\Omega _2(t)} u\, dxdy$ is conserved for any harmonic function $u(x,y)$ defined in $\Omega$. A special choice of $u(x,y)\equiv 1$ for $\chi_2=0, 1$, corresponds to the volume conservation - in that case, the integral on the right hand side is zero.

Remark that in the Saffman-Taylor formulation of the problem -  where a viscous fluid occupying the gap between two plates is being displaced by a less viscous fluid,  which is forced into the gap - unstable fingers are being formed. Similarly, a basic instability - a version of the Saffman-Taylor instability - was identified in \cite{tian} when a viscous circular bubble was surrounded by the air and the upper plate was lifting.

Unstable fingers are subject to tip splitting and exhibit singularities in a finite time.
In the present paper we did not consider neither  formation of singularities, nor the ways of achieving a regularization. The aim of this study was, in contrast, to avoid formation of singularities by means of a special choice of sinks and sources. 
Note that linear stability results for the interior problem \cite{tian}
indicate that a circular bubble is stable
when the plate is moving down. The latter, together with the  stability results for the Saffman-Taylor formulation in a radial flow geometry \cite{miranda}, suggest 
that  the circular interface for the problem in question is expected to be linearly stable in two situations: (i) when a more viscous fluid occupies the interior domain and the upper plate is moving down or (ii) when a less viscous fluid is surrounded by a more viscous fluid and the upper plate is moving up.




\providecommand{\bysame}{\leavevmode\hbox to3em{\hrulefill}\thinspace}

\end{document}